\documentclass[10pt]{article}
\usepackage{amssymb}
\usepackage{amsmath}
\usepackage[dvips]{graphics}
\usepackage{epsfig}

\begin{document}
\title{A Probabilistic Model for
Forest Fires}

\author{\textbf{Vassilis G. Papanicolaou}
\\\\
Department of Mathematics
\\
National Technical University of Athens,
\\
Zografou Campus, 157 80, Athens, GREECE
\\
{\tt papanico@math.ntua.gr}}
\maketitle

\begin{abstract}
We propose a discrete two-dimensional mathematical model for forest fires and we 
derive certain results describing its limiting behavior. We also pose a relevant
open question.
\end{abstract}

\textbf{Keywords.} Random field; limiting behavior.
\\\\
\textbf{2010 AMS Mathematics Classification.} 60G60; 92F05.

\section{Introduction}
The forest is represented by the set
\begin{equation*}
\mathbb{N}^2 = \{(j, k) \,:\, j, k \in \mathbb{N}\},
\end{equation*}
where $\mathbb{N} = \{0, 1, 2, \ldots\}$ is the set of nonnegative integers.
Each point $(j, k) \in 
\mathbb{N}^2$ is associated to a tree. The \emph{status} of the tree at $(j, k)$ is denoted by 
$X(j, k)$ (later we will use additional notations for $X(j, k)$), and it takes two 
values, $0$ or $1$, so that $X(j, k) = 1$ means that the tree at $(j, k)$ is burnt, while
$X(j, k) = 0$ means that the tree at $(j, k)$ is not burnt.

The key feature of the model is the following: The status of the tree at $(j, k)$ is affected only by 
the status of the neighboring trees at $(j-1, k)$ and $(j, k-1)$ so that

(i) if the trees at $(j-1, k)$ and $(j, k-1)$ are not burnt, then tree at $(j, k)$ is not burnt;

(ii) if the tree at $(j-1, k)$ is burnt and the tree at $(j, k-1)$ is not burnt,
then there is a probability $\alpha$ that the tree at $(j, k)$ is burnt;

(iii) if the tree at $(j-1, k)$ is not burnt and the tree at $(j, k-1)$ is burnt,
then there is a probability $\beta$ that the tree at $(j, k)$ is burnt;

(iv) if both trees at $(j-1, k)$ and $(j, k-1)$ are burnt,
then there is a probability $\gamma$ that the tree at $(j, k)$ is burnt.

To avoid paradoxical situations we need to assume that $\gamma \geq \max\{\alpha, \beta\}$, since one 
expects that in the
case where both trees at $(j-1, k)$ and $(j, k-1)$ are burnt there are more chances that
the tree at $(j, k)$ is burnt than in the case where only one of the trees at $(j-1, k)$ 
and $(j, k-1)$ is burnt.

The probabilities $\alpha, \beta, \gamma$ are related to the wind speed and wind direction.


We assume that the fire starts at the point $(0, 0)$ so that $X(0, 0) = 1$.

The complete description of the model is given in the next section.

The main results of the paper concern the sequence of random variables
\begin{equation*}
Y_n = X(0, n) + X(1, n-1) + \cdots + X(n, 0)
\end{equation*}
since this sequence contains significant information regarding the fate of the forest.

In the case $\alpha + \beta < 1$ and  $\gamma < 1$ we show that
$\sum_{n=0}^{\infty} Y_n < \infty$ almost surely,
which implies that the forest is not considerably damaged.

The case $\alpha + \beta = \gamma = 1$, with $\alpha \beta > 0$, is a borderline case. Here we show that
$Y_n \to 0$ almost surely, but $\mathbb{E}[Y_n] = 1$ for all $n \geq 0$. This gives a somehow ``mixed"
information regarding the fate of the forest.

The case $\alpha + \beta > 1$ is the most challenging. Here, with the additional hypothesis that
$\gamma = 1$ we show that
\begin{equation*}
\mathbb{P}\left\{Y_n = 0 \right\} \leq \frac{n}{n+1} (2 - \alpha - \beta)
\qquad
\text{and}
\qquad
\mathbb{E}\left[Y_n \right] \geq 1 + (\alpha + \beta - 1)n,
\qquad
n \geq 0.
\end{equation*}
These estimates indicate that the fire is quite catastrophic.

In the opposite direction, we have the bounds
\begin{equation*}
\lim_n\mathbb{P}\left\{Y_n = 0 \right\} \geq \frac{(1 - \alpha)(1 - \beta)}{\alpha \beta}
\qquad
\text{and}
\qquad
\limsup_n \frac{\mathbb{E}\left[Y_n \right]}{n} \leq 
\frac{\alpha + \beta - 1}{\alpha \beta},
\end{equation*}
which are, clearly, also valid in the case where $\gamma < 1$ (as long as $\alpha + \beta > 1$).

Finally, in the case $\alpha + \beta > 1$ an interesting open question is whether $Y_n/n$ converges, at least in distribution.

\section{The law of the random field $X(j, k)$}
In the sequel we will use the notations
\begin{equation}
X_j^{j+k} = X_{jk} = X(j, k),
\qquad
\text{so that}
\qquad
X_j^n = X(j, n-j),
\label{B0}
\end{equation}
e.g., $X_2^5 = X_{23} = X(2, 3)$. Thus, for convenience, sometimes $X(j, k)$ will be denoted by 
$X_{jk}$ and sometimes by $X_j^n$, where $n = j + k$.

The ``boundary conditions" for $X(j, k)$ are
\begin{equation}
X(0, 0) = 1
\qquad
\text{and}
\qquad
X(j, -1) = X(-1, k) = 0
\qquad
\text{for all } j,k \in \mathbb{N}
\label{A4}
\end{equation}
(we can further assume that $X(j, k) = 0$ whenever $j < 0$ or $k < 0$).

We will now define the random field $X(j, k)$, $j,k \in \mathbb{N}$, inductively.

Notice that $X(j, k)$ is a Bernoulli random variable, i.e. it takes only the values
$0$ and $1$.

Start with the boundary conditions \eqref{A4}. Next, suppose that for some $n \geq 0$ we have defined 
the field $X_j^{j+k} = X_{jk}$ for $0 \leq j+k \leq n$. We will show how to define (the law of)
$X_j^{j+k}$ for $j+k = n+1$.

First we introduce the $\sigma$-algebras
\begin{equation}
\mathcal{F}_n = \sigma\left(X_{jk},\ 0 \leq j+k \leq n\right),
\qquad
\mathcal{G}_n = \sigma\left(X_{jk},\  j+k = n\right),
\qquad
n \geq 0
\label{B1}
\end{equation}
(of course, in view of \eqref{A4}, we have that
$\mathcal{F}_0 = \mathcal{G}_0 = \{\varnothing, \Omega\}$).

For each (fixed) $j \in \{0,1, \ldots, n+1\}$ we define the law of $X_j^{n+1}$ by setting
\begin{equation}
\left.
  \begin{array}{l}
\mathbb{P}\big\{X_j^{n+1} = 1\,\big|\, X_{j-1}^n = X_j^n = 0\big\} = 0,
\\
\mathbb{P}\big\{X_j^{n+1} = 1\,\big|\, X_{j-1}^n = 1, \ X_j^n = 0\big\} = \alpha,
\\
\mathbb{P}\big\{X_j^{n+1} = 1\,\big|\, X_{j-1}^n = 0, \ X_j^n = 1\big\} = \beta,
\\
\mathbb{P}\big\{X_j^{n+1} = 1\,\big|\, X_{j-1}^n = X_j^n = 1\big\} = \gamma,
  \end{array}
\right\}
\label{A1}
\end{equation}
where
\begin{equation}
0 \leq \alpha, \beta \leq \gamma \leq 1
\label{A2}
\end{equation}
($\alpha$, $\beta$, and $\gamma$ are given). The set of equations in \eqref{A1} presents the main feature of the 
model, namely that the status of the tree at $(j, k)$ is affected only by the status of the 
neighboring trees at $(j-1, k)$ and $(j, k-1)$.

Now, the joint distribution of $X_j^{n+1}$, $0 \leq j \leq n+1$, is defined by setting

\begin{align}
\mathbb{P}\big\{X_j^{n+1} = \epsilon_j, \ 0 \leq j \leq n+1 \, \big| \,
\mathcal{F}_n\big\}
&= \mathbb{P}\big\{X_j^{n+1} = \epsilon_j, \ 0 \leq j \leq n+1 \, \big| \,  
\mathcal{G}_n\big\}
\nonumber
\\
&= \prod_{j=0}^{n+1} \mathbb{P}\big\{X_j^{n+1} = \epsilon_j
\, \big| \, \mathcal{G}_n\big\}
\nonumber
\\
&= \prod_{j=0}^{n+1} \mathbb{P}\big\{X_j^{n+1} = \epsilon_j \, \big| \, X_{j-1}^n, X_j^n\big\},
\label{B2}
\end{align}
where each $\epsilon_j$, $0 \leq j \leq n+1$, can be either $0$ or $1$ (recall that $X(-1, n)$ and 
$X(n, -1)$ are deterministic, i.e. they have given nonrandom values).

The second equality in \eqref{B2} states that the random variables $X_j^{n+1}$, $0 \leq j \leq n+1$, are conditionally independent given $\mathcal{G}_n$. Actually, since \eqref{B2} implies that 
\begin{equation*}
\mathbb{P}\big\{X_j^{n+1} = \epsilon_j \, \big| \, \mathcal{G}_n\big\}
= \mathbb{P}\big\{X_j^{n+1} = \epsilon_j \, \big| \, \mathcal{F}_n\big\},
\end{equation*}
it follows that $X_j^{n+1}$, $0 \leq j \leq n+1$, are also conditionally independent given
$\mathcal{F}_n$ (for instance, since $\mathcal{F}_0$ is the trivial $\sigma$-algebra, the random variables
$X_0^1$ and $X_1^1$ are independent).

Let us also notice that the defining equations in \eqref{B2} have a Markovian flavor \cite{K}.

One consequence of the definition of the random field $X_j^n$ is that there is a subset $\Omega_0$ of 
$\Omega$ with $\mathbb{P}(\Omega_0) = 1$ such that
\begin{equation}
\left\{X_j^n = 0,\ 0 \leq j \leq n \right\} \cap \Omega_0
\subset \left\{X_j^{n+1} = 0,\ 0 \leq j \leq n+1\right\} \cap \Omega_0,
\qquad
n \geq 0.
\label{BA2}
\end{equation}

\medskip

\textbf{Remark 1.} From the definition of the random field $X(j, k)$ it follows that if
\begin{equation*}
\tilde{\alpha} \geq \alpha, 
\qquad
\tilde{\beta} \geq \beta,
\qquad
\tilde{\gamma} \geq \gamma
\end{equation*}
and 
$\tilde{X}(j, k)$ is the random field associated to
$\tilde{\alpha}, \tilde{\beta}, \tilde{\gamma}$, while
$X(j, k)$ is the random field associated to $\alpha, \beta, \gamma$, then
\begin{equation}
\mathbb{P}\left\{\tilde{X}(j, k) = 1\right\} \geq \mathbb{P}\left\{X(j, k) = 1\right\}
\qquad
\text{for every }
j,k \geq 0,
\label{B2bb}
\end{equation}
i.e. $\tilde{X}(j, k)$ is stochastically larger than $X(j, k)$ for every $j,k \geq 0$.
\hfill $\diamondsuit$

\medskip

The symbol $\diamondsuit$ indicates the end of a remark or an example.

\medskip

\textbf{Remark 2.} From formulas \eqref{A4} and \eqref{A1} it follows that
if $\alpha = \beta = 1$, then $X(j, k) = 1$ a.s. for all $(j,k) \in \mathbb{N}^2$,
i.e. all trees of the forest are burnt, while if $\alpha = \beta = 0$,
then $X(j, k) = 0$ a.s.
for all $(j,k) \in \mathbb{N}^2 \smallsetminus \{(0,0)\}$, i.e. the only burnt tree is 
the one at $(0,0)$. 

If $\alpha = 0$, then \eqref{A4} and \eqref{A1} imply that the only trees that are possibly burnt
are located at the points $(0, k)$, $k \in \mathbb{N}$. Likewise, if $\beta = 0$, then \eqref{A4} 
and \eqref{A1} imply that the only trees that are possibly burnt
are located at the points $(j, 0)$, $j \in \mathbb{N}$. Therefore, if $\alpha \beta = 0$, our
model is one-dimensional, and, consequently, very simple. We analyze this model in the Appendix.
\hfill $\diamondsuit$

\medskip

For the rest of the paper we will always assume that
\begin{equation}
\alpha > 0
\qquad
\text{and}
\qquad
\beta > 0.
\label{B2bbb}
\end{equation}

Next, we present some immediate consequences of the defining formulas \eqref{A1} and 
\eqref{B2}. But, first, let us recall that for a Bernoulli random variable $X$ we have that
\begin{equation}
{\bf 1}_{\{X = 1\}} = X,
\qquad
{\bf 1}_{\{X = 0\}} = 1 - X,
\quad
\text{and}
\quad
\mathbb{E}\left[z^X\right] = 1 + \mathbb{E}[X] (z-1).
\label{B2b}
\end{equation}

From \eqref{A1}, \eqref{B2}, and \eqref{B2b} it follows that
\begin{align}
\mathbb{P}\left\{X_j^n = 1 \, \big| \, \mathcal{F}_{n-1}\right\}
&=\mathbb{E}\left[X_j^n \, \big| \, \mathcal{F}_{n-1}\right]
\nonumber
\\
&=\mathbb{E}\left[X_j^n \, \big| \, X_{j-1}^{n-1}, X_j^{n-1}\right]
\nonumber
\\
&=\alpha X_{j-1}^{n-1} \left(1 - X_j^{n-1}\right) 
+ \beta X_j^{n-1} \left(1 - X_{j-1}^{n-1}\right)
+ \gamma X_{j-1}^{n-1} X_j^{n-1}
\nonumber
\\
&=\alpha X_{j-1}^{n-1} + \beta X_j^{n-1} -(\alpha + \beta - \gamma) X_{j-1}^{n-1} X_j^{n-1}.
\label{B2c}
\end{align}
Also, from \eqref{B2b} and \eqref{B2c} we get
\begin{align}
\mathbb{E}\left[z^{X_j^n} \, \big| \, \mathcal{F}_{n-1}\right]
&= 1 + \mathbb{E}\left[X_j^n \, \big| \, \mathcal{F}_{n-1}\right] (z-1)
\nonumber
\\
&= 1 + \left[\alpha X_{j-1}^{n-1} + \beta X_j^{n-1}
-(\alpha + \beta - \gamma) X_{j-1}^{n-1} X_j^{n-1}\right] (z-1).
\label{B2d}
\end{align}
Finally, from \eqref{B2c} and the fact that, for $j \ne k$, the random variables
$X_j^n$ and $X_k^n$ are conditionally independent given $\mathcal{F}_{n-1}$ we get
\begin{align}
\mathbb{E}\left[X_j^n X_k^n \, \big| \, \mathcal{F}_{n-1}\right]
&=\mathbb{E}\left[X_j^n \, \big| \, \mathcal{F}_{n-1}\right]
\mathbb{E}\left[ X_k^n \, \big| \, \mathcal{F}_{n-1}\right]
\nonumber
\\
&=\left[\,\alpha X_{j-1}^{n-1} + \beta X_j^{n-1} 
-(\alpha + \beta - \gamma) X_{j-1}^{n-1} X_j^{n-1}\right]
\nonumber
\\
&\quad \times \left[\,\alpha X_{k-1}^{n-1} + \beta X_k^{n-1} 
-(\alpha + \beta - \gamma) X_{k-1}^{n-1} X_k^{n-1}\right].
\label{B2e}
\end{align}

At the end of this section we give few examples which indicate how to calculate certain
probabilities related to the field $X(j, k)$.

\medskip

\textbf{Example 1.} For a given $n \geq 0$ let us show how to compute the probability
\begin{equation}
\mathbb{P}\left\{X_j^n = \epsilon_j, \ 0 \leq j \leq n \right\},
\label{B3}
\end{equation}
where the value of each $\epsilon_j$, $0 \leq j \leq n$, is given
(it is either $0$ or $1$).

For $n=0$ we have that $X_0^0 = 1$ is deterministic.

For $n=1$ the probability in \eqref{B3} becomes, in view of \eqref{B2} and the 
fact that $\mathcal{G}_0$ is the trivial $\sigma$-algebra
\begin{equation}
\mathbb{P}\left\{X_0^1 = \epsilon_0, \ X_1^1 = \epsilon_1 \right\}
= \mathbb{P}\left\{X_0^1 = \epsilon_0 \right\}
\mathbb{P}\left\{X_1^1 = \epsilon_1 \right\}
\label{B3a}
\end{equation}
and these probabilities can be computed from the boundary conditions 
together with the equations of \eqref{A1}.

Suppose that for $n \leq m$ we can compute the probability in \eqref{B3} by using
\eqref{B2} and \eqref{A1}
(and the boundary condition). We will, then, show how to calculate the probability 
in \eqref{B3} for $n = m+1$. We have
\begin{align}
&\mathbb{P}\left\{X_j^{m+1} = \epsilon_j, \ 0 \leq j \leq m+1 \right\}
\nonumber
\\
\qquad \qquad \  &= \sum_{\tilde{\epsilon}_k = 0 \text{ or } 1 \atop 0 \leq k \leq m}
\left[\,\mathbb{P}\left\{X_j^{m+1} = \epsilon_j, \ 0 \leq j \leq m+1 \, \big| \,
X_k^m = \tilde{\epsilon}_k, \ 0 \leq k \leq m\right\}\right.
\nonumber
\\
&\qquad \qquad \ \ \left. \times\, \mathbb{P}\left\{X_k^m = \tilde{\epsilon}_k, \ 0 \leq k \leq m \right\}\right],
\label{B4}
\end{align}
where the sum is taken over all the $2^{m+1}$ different choices of
$\tilde{\epsilon}_k$, $0 \leq k \leq m$.

Each term of the sum in \eqref{B4} is a product of a conditional probability and a
probability. The latter can be computed by the induction hypothesis. As for
the conditional probability, it can be computed too by using \eqref{B2} (actually, in view of
\eqref{B2}, the events $\{X_j^{m+1} = \epsilon_j\}$, $0 \leq j \leq m+1$,
are conditionally independent given the event $\{X_k^m = \tilde{\epsilon}_k, \ 0 \leq k \leq m\}$).
This finishes the induction which shows how all probabilities of the
form of \eqref{B3} can be computed.
\hfill $\diamondsuit$

\medskip

\textbf{Example 2.} Let us compute the probability
\begin{equation}
\mathbb{P}\big\{X_{12} = 0, \  X_{11} = 1 \big\}.
\label{B5}
\end{equation}
We have
\begin{equation*}
\mathbb{P}\big\{X_{12} = 0, \  X_{11} = 1 \big\}
= \sum_{\epsilon = 0, 1}
\mathbb{P}\big\{X_{12} = 0 \, \big| \, X_{02} = \epsilon, \  X_{11} = 1 \big\}
\mathbb{P}\big\{X_{02} = \epsilon, \  X_{11} = 1 \big\},
\end{equation*}
where all quantities in the right-hand side can be computed. Indeed, \eqref{A1} 
yields
\begin{equation*}
\mathbb{P}\big\{X_{12} = 0 \, \big| \, X_{02} = 0, \  X_{11} = 1 \big\} = 1 - \beta
\end{equation*}
and
\begin{equation*}
\mathbb{P}\big\{X_{12} = 0 \, \big| \, X_{02} = 1, \  X_{11} = 1 \big\} = 1 - \gamma,
\end{equation*}
while the probability $\mathbb{P}\{X_{02} = \epsilon, \  X_{11} = 1 \}$
can be computed with the help of Example 1.
\hfill $\diamondsuit$

\medskip

\textbf{Example 3.} Let us also compute the probability
\begin{equation}
\mathbb{P}\big\{X_{22} = 1, \ X_{12} = 0, \  X_{11} = 1 \big\}.
\label{B6}
\end{equation}
We have
\begin{align}
&\mathbb{P}\big\{X_{22} = 1, \ X_{12} = 0, \  X_{11} = 1 \big\}
\nonumber
\\
= &\sum_{\epsilon = 0, 1}
\mathbb{P}
\big\{X_{22} = 1 \, \big| \, X_{12} = 0, \ X_{21} = \epsilon, \  X_{11} = 1 \big\}
\mathbb{P}\big\{X_{12} = 0, \ X_{21} = \epsilon, \  X_{11} = 1 \big\},
\nonumber
\\
= &\sum_{\epsilon = 0, 1}
\mathbb{P}
\big\{X_{22} = 1 \, \big| \, X_{12} = 0, \ X_{21} = \epsilon \big\}
\mathbb{P}\big\{X_{12} = 0, \ X_{21} = \epsilon, \  X_{11} = 1 \big\}
\nonumber
\\
= &\,\mathbb{P}
\big\{X_{22} = 1 \, \big| \, X_{12} = 0, \ X_{21} = 1 \big\}
\mathbb{P}\big\{X_{12} = 0, \ X_{21} = 1, \  X_{11} = 1 \big\}
\nonumber
\\
= &\,\beta \,\mathbb{P}\big\{X_{12} = 0, \ X_{21} = 1, \  X_{11} = 1 \big\},
\nonumber
\end{align}
where the second equality is a consequence of \eqref{B2}, while the last two 
equalities follow from \eqref{A1}.

Finally,
\begin{align}
&\mathbb{P}\big\{X_{12} = 0, \ X_{21} = 1, \ X_{11} = 1 \big\}
\nonumber
\\
&= \sum_{\epsilon_0, \epsilon_2 = 0, 1}
\left[\,\mathbb{P}\big\{X_{12} = 0, \ X_{21} = 1 \,\big| \, 
X_{02} = \epsilon_0, \ X_{11} = 1,\ X_{20} = \epsilon_2 \big\} \right.
\nonumber
\\
&\qquad\qquad\ \left.\times \,
\mathbb{P}\big\{X_{02} = \epsilon_0, \ X_{11} = 1,\ X_{20} = \epsilon_2 \big\} \right]
\nonumber
\\
&= \sum_{\epsilon_0, \epsilon_2 = 0, 1}
\left[\,\mathbb{P}\big\{X_{12} = 0 \,\big| \, 
X_{02} = \epsilon_0, \ X_{11} = 1 \big\}
\mathbb{P}\big\{X_{21} = 1 \,\big| \, 
X_{11} = 1, \ X_{20} = \epsilon_2 \big\} \right.
\nonumber
\\
&\qquad\qquad\ \left.\times \,
\mathbb{P}\big\{X_{02} = \epsilon_0, \ X_{11} = 1, X_{20} = \epsilon_2 \big\} \right]
\nonumber
\end{align}
(the last equality follows from \eqref{B2}, namely the conditional independence given
$\mathcal{G}_2$), where all probabilities can be computed with the help of \eqref{A1} and Example 1.
\hfill $\diamondsuit$

\medskip

\textbf{Example 4.} Here we compute the conditional probabilities
\begin{equation}
w_k = \mathbb{P}\left\{X_1^k = 0 \,\big| \, X_0^1 = X_0^2 = \cdots = X_0^k = 1\right\},
\qquad
k \geq 1.
\label{B7}
\end{equation}
For $k=1$ we have (since $X_0^1$ and $X_1^1$ are independent)
\begin{equation}
w_1 = \mathbb{P}\left\{X_1^1 = 0 \,\big| \, X_0^1= 1\right\} = \mathbb{P}\left\{X_1^1 = 0\right\}
= 1 - \alpha.
\label{B8}
\end{equation}

For $k \geq 2$ the conditional probability of \eqref{B7} can be expressed as
\begin{equation}
w_k 
 = \frac{\mathbb{P}\left\{X_0^k = 1, X_1^k = 0 \,\big| \, X_0^1 = X_0^2 = \cdots = X_0^{k-1} = 1\right\}}{\mathbb{P}\left\{X_0^k = 1\,\big| \, X_0^1 = X_0^2 = \cdots = X_0^{k-1} = 1\right\}},
\label{B9}
\end{equation}
where from the defining properties of the random field $X_j^n$ we have that
\begin{equation}
\mathbb{P}\left\{X_0^k = 1\,\big| \, X_0^1 = X_0^2 = \cdots = X_0^{k-1} = 1\right\}
= \mathbb{P}\left\{X_0^k = 1\,\big| \, X_0^{k-1} = 1\right\} = \beta.
\label{B10}
\end{equation}
hence \eqref{B9} becomes
\begin{equation}
w_k = \beta^{-1}
\mathbb{P}\left\{X_0^k = 1, X_1^k = 0 \,\big| \, A_{k-1}\right\},
\label{B11}
\end{equation}
where for typographical convenience we have set
\begin{equation}
A_k = \left\{X_0^1 = X_0^2 = \cdots = X_0^k = 1\right\},
\qquad
k \geq 1.
\label{B11a}
\end{equation}

Now,
\begin{align}
&\mathbb{P}\left\{X_0^k = 1, X_1^k = 0 \,\big| \, A_{k-1}\right\}
\nonumber
\\
&= \mathbb{P}\left\{X_0^k = 1, X_1^k = 0 \,\big| \, A_{k-1}, X_1^{k-1} = 0\right\}
\mathbb{P}\left\{X_1^{k-1} = 0 \,\big| \, A_{k-1}\right\}
\nonumber
\\
&\quad + \mathbb{P}\left\{X_0^k = 1, X_1^k = 0 \,\big| \, A_{k-1}, X_1^{k-1} = 1\right\}
\mathbb{P}\left\{X_1^{k-1} = 1 \,\big| \, A_{k-1}\right\}.
\label{B12}
\end{align}
In view of \eqref{B2}, \eqref{B7}, and \eqref{B11a}, formula  \eqref{B12} becomes
\begin{align}
&\mathbb{P}\left\{X_0^k = 1, X_1^k = 0 \,\big| \, A_{k-1}\right\}
\nonumber
\\
&= \mathbb{P}\left\{X_0^k = 1, X_1^k = 0 \,\big| \, X_0^{k-1} = 1, X_1^{k-1} = 0\right\}
w_{k-1}
\nonumber
\\
&\quad + \mathbb{P}\left\{X_0^k = 1, X_1^k = 0 \,\big| \, X_0^{k-1} = 1, X_1^{k-1} = 1\right\}
(1 - w_{k-1})
\nonumber
\\
&= \mathbb{P}\left\{X_0^k = 1 \,\big| \, X_0^{k-1} = 1\right\}
\mathbb{P}\left\{X_1^k = 0 \,\big| \, X_0^{k-1} = 1, X_1^{k-1} = 0\right\}
w_{k-1}
\nonumber
\\
&\quad + \mathbb{P}\left\{X_0^k = 1 \,\big| \, X_0^{k-1} = 1\right\}
\mathbb{P}\left\{X_1^k = 0 \,\big| \, X_0^{k-1} = 1, X_1^{k-1} = 1\right\}
(1 - w_{k-1}).
\label{B13}
\end{align}
Thus, by using \eqref{A1} in \eqref{B13} we obtain
\begin{equation}
\mathbb{P}\left\{X_0^k = 1, X_1^k = 0 \,\big| \, A_{k-1}\right\}
= \beta (1 - \alpha)w_{k-1} + \beta (1 - \gamma)(1 - w_{k-1}),
\label{B14}
\end{equation}
and, finally, in view of \eqref{B14}, formula \eqref{B11} yields
\begin{equation}
w_k = (1 - \alpha)w_{k-1} + (1 - \gamma)(1 - w_{k-1}).
\label{B15}
\end{equation}
Formula \eqref{B15} is a linear first-order difference equation with constant coefficients, whose
initial condition is \eqref{B8}. Therefore, it is easy to get that
\begin{equation}
w_k = \frac{\alpha}{1 - (\gamma - \alpha)} (\gamma - \alpha)^k
+ \frac{1 - \gamma}{1 - (\gamma - \alpha)},
\qquad
k \geq 1.
\label{B16}
\end{equation}
\hfill $\diamondsuit$

%

\subsection{An equivalent Markovian model}
The above forest fire model can be viewed as a Markov chain.

Let us start with a definition.

\medskip

\textbf{Definition 1.} A \emph{binary sequence} is a sequence
$b = \{b_j\}_{j=0}^{\infty}$ such that $b_j = 0 $ or $1$ for every $j \geq 0$.

The \emph{valence} $\boldsymbol{v}(b)$ of a binary sequence
$b = \{b_j\}_{j=0}^{\infty}$ is the 
number of its elements which are equal to $1$, namely
\begin{equation}
\boldsymbol{v}(b) = \sum_{j=0}^{\infty} b_j.
\label{N0a}
\end{equation}
The set of all binary sequences of finite valence is denoted by $\mathcal{S}$. Thus
\begin{equation}
\mathcal{S}
= \left\{b \,:\, b \text{ is a binary sequence with } \boldsymbol{v}(b)
< \infty\right\}.
\label{N0b}
\end{equation}

\medskip

Next, we define a Markov chain $\hat{X}(n)$, $n \in \mathbb{N}$,
whose state space is the set $\mathcal{S}$ of \eqref{N0b}. Thus, for each
$n \in \mathbb{N}$ we have that
\begin{equation}
\hat{X}(n) = \left\{\hat{X}_j(n)\right\}_{j=0}^{\infty}
\quad \text{with} \quad
\hat{X}_j(n) = 0 \text{ or } 1 
\quad \text{and} \quad
\sum_{j=0}^{\infty} \hat{X}_j(n) < \infty.
\label{N1}
\end{equation}
The transition probabilities of $\hat{X}(n)$ are given by adapting
\eqref{A1} and \eqref{B2}. Thus, for every $j,n \in \mathbb{N}$ we have that
\begin{equation}
\mathbb{P}\left\{\hat{X}_j(n+1) \,\big|\, \hat{X}(n)\right\}
= \mathbb{P}\left\{\hat{X}_j(n+1) \,\big|\, \hat{X}_{j-1}(n),\ \hat{X}_j(n)\right\}
\label{N2}
\end{equation}
(with the convention that $\hat{X}_{-1}(n) = 0$ for all $n \in \mathbb{N}$) 
and (compare with \eqref{A1})
\begin{equation}
\left.
  \begin{array}{l}
\mathbb{P}\left\{\hat{X}_j(n+1) = 1\,\big|\, \hat{X}_{j-1}(n) = \hat{X}_j(n) = 0\right\} = 0,
\\
\mathbb{P}\left\{\hat{X}_j(n+1) = 1\,\big|\, \hat{X}_{j-1}(n) = 1, \ \hat{X}_j(n) = 0\right\} = \alpha,
\\
\mathbb{P}\left\{\hat{X}_j(n+1) = 1\,\big|\, \hat{X}_{j-1}(n) = 0, \ \hat{X}_j(n) = 1\right\} = \beta,
\\
\mathbb{P}\left\{\hat{X}_j(n+1) = 1\,\big|\, \hat{X}_{j-1}(n) = \hat{X}_j(n) = 1\right\} = \gamma.
  \end{array}
\right\}
\label{N3}
\end{equation}
As before, $\alpha$, $\beta$, and $\gamma$ are given constants satisfying \eqref{A2}.
To complete the description of the law of $\hat{X}(n)$ we pose that the (Bernoulli) 
random variables $\hat{X}_0(n+1), \hat{X}_1(n+1), \hat{X}_2(n+1), \ldots$ are conditionally independent
given $\hat{X}(n)$, for every $n \in \mathbb{N}$
(this assumption, together with \eqref{N2}, is equivalent to \eqref{B2}).

The boundary conditions \eqref{A4} correspond to the initial condition
\begin{equation}
\hat{X}_0(0) = 1
\qquad
\text{and}
\qquad
\hat{X}_j(0) = 0
\quad
\text{for all } j \geq 1.
\label{N4}
\end{equation}

Finally, let us mention that the valences
\begin{equation}
\hat{Y}(n) = \boldsymbol{v}\left(\hat{X}(n)\right)
= \sum_{j=0}^{\infty} \hat{X}_j(n),
\qquad
n \in \mathbb{N},
\label{N5}
\end{equation}
correspond to the variables $Y_n$, $n \in \mathbb{N}$, introduced in the next 
section.

\section{The random variables $Y_n$, $n \geq 0$}
It is somehow natural to introduce the random variables
\begin{equation}
Y_n = \sum_{j=0}^n X_j^n,
\qquad
n \geq 0.
\label{C1}
\end{equation}
Actually, regarding the forest fire, the limiting behavior, as $n \to \infty$, of $Y_n$ 
plays a dominant role in the fate of the forest.

Notice that
\begin{equation}
Y_0 = 1
\qquad
\text{and}
\qquad
Y_n \in \{0, 1, \ldots, n+1\}.
\label{C1a}
\end{equation}
Also, by the defining properties of the random field $X_j^n$ (recall \eqref{BA2}) there is a subset
$\Omega_0$ of $\Omega$ with $\mathbb{P}(\Omega_0) = 1$ such that
\begin{equation}
\left\{Y_n = 0 \right\} \cap \Omega_0 \subset \left\{Y_{n+1} = 0 \right\} \cap \Omega_0,
\qquad
n \geq 0,
\label{C1b}
\end{equation}
and, more generally,
\begin{equation}
\left\{Y_n \leq m \right\} \cap \Omega_0 \subset \left\{Y_{n+1}  \leq 2m \right\} \cap \Omega_0,
\qquad
n \geq 0,
\label{C1c}
\end{equation}

We have
\begin{equation}
\mathbb{E}\left[Y_n \, \big| \, \mathcal{F}_{n-1}\right]
= \sum_{j=0}^n \mathbb{E}\left[X_j^n \, \big| \, \mathcal{F}_{n-1}\right],
\qquad
n \geq 1,
\label{C2a}
\end{equation}
hence, in view of \eqref{B2c},
\begin{equation*}
\mathbb{E}\left[Y_n \, \big| \, \mathcal{F}_{n-1}\right]
=\alpha \sum_{j=0}^n X_{j-1}^{n-1} + \beta \sum_{j=0}^n X_j^{n-1}
-(\alpha + \beta - \gamma) \sum_{j=0}^n X_{j-1}^{n-1} X_j^{n-1},
\nonumber
\end{equation*}
from which, and the boundary conditions \eqref{A4}, it follows that
\begin{equation}
\mathbb{E}\left[Y_n \, \big| \, \mathcal{F}_{n-1}\right]
=(\alpha + \beta) Y_{n-1}
-(\alpha + \beta - \gamma) \sum_{j=1}^{n-1} X_{j-1}^{n-1} X_j^{n-1},
\qquad
n \geq 1.
\label{C2}
\end{equation}

\section{The case $\alpha + \beta = \gamma$}
Formula \eqref{C2} suggests that the case
\begin{equation}
\alpha + \beta = \gamma
\label{F1}
\end{equation}
will be easier to analyze. In this section it is always assumed that \eqref{F1} is satisfied.

Under \eqref{F1}, formula \eqref{B2c} becomes
\begin{equation}
\mathbb{P}\left\{X_j^n = 1 \, \big| \, \mathcal{F}_{n-1}\right\}
= \mathbb{E}\left[X_j^n \, \big| \, \mathcal{F}_{n-1}\right]
=\alpha X_{j-1}^{n-1} + \beta X_j^{n-1}.
\label{F2}
\end{equation}
Since $\mathcal{F}_{n-2} \subset \mathcal{F}_{n-1}$, by conditioning on $\mathcal{F}_{n-2}$
formula \eqref{F2} yields
\begin{align}
\mathbb{E}\left[X_j^n \, \big| \, \mathcal{F}_{n-2}\right]
&=\alpha \mathbb{E}\left[X_{j-1}^{n-1} \, \big| \, \mathcal{F}_{n-2}\right]
+ \beta \mathbb{E}\left[X_j^{n-1} \, \big| \, \mathcal{F}_{n-2}\right]
\nonumber
\\
&=\alpha \left[\alpha X_{j-2}^{n-2} + \beta X_{j-1}^{n-2}\right]
+ \beta \left[\alpha X_{j-1}^{n-2} + \beta X_j^{n-2}\right]
\\
&= \alpha^2 X_{j-2}^{n-2} + 2 \alpha \beta X_{j-1}^{n-2} + \beta^2 X_j^{n-2}
\label{F3}
\end{align}
and by straightforward induction
\begin{equation}
\mathbb{P}\left\{X_j^n = 1 \, \big| \, \mathcal{F}_{n-m}\right\}
= \mathbb{E}\left[X_j^n \, \big| \, \mathcal{F}_{n-m}\right]
= \sum_{k=0}^m \left({m} \atop {k}\right) \alpha^k \beta^{m-k} X_{j-k}^{n-m},
\quad
0 \leq m \leq n
\label{F4}
\end{equation}
(in the case $m=0$, formula \eqref{F4} is trivially true since
$X_j^n$ is $\mathcal{F}_n$-measurable).

In particular, for $m=n$, in view of \eqref{A4} and the fact that
$\mathcal{F}_0 = \{\varnothing, \Omega\}$, formula \eqref{F4} yields
\begin{equation}
\mathbb{P}\left\{X_j^n = 1\right\} = \mathbb{E}\left[X_j^n\right]
= \sum_{k=0}^n \left({n} \atop {k}\right) \alpha^k \beta^{n-k} X_{j-k}^0
= \left({n} \atop {j}\right) \alpha^j \beta^{n-j}.
\label{F5}
\end{equation}
Formula \eqref{F5} remains valid even in the case where $j$ is not between $0$ and $n$,
since if $j < 0$ or $j > n$, then $X_j^n = X_{j(n-j)} = 0$ and
at the same time the binomial coefficient $\left({n} \atop {j}\right)$ vanishes.

Now, regardind $Y_n$ we notice that, under \eqref{F1}, formula \eqref{C2} becomes
\begin{equation}
\mathbb{E}\left[Y_n \, \big| \, \mathcal{F}_{n-1}\right]
=(\alpha + \beta) Y_{n-1}.
\label{F6}
\end{equation}
It follows that if we set
\begin{equation}
M_n = (\alpha + \beta)^{-n} Y_n,
\qquad
n \geq 0,
\label{F7}
\end{equation}
then
\begin{equation}
\mathbb{E}\left[M_n \, \big| \, \mathcal{F}_{n-1}\right]
= (\alpha + \beta)^{-n} \mathbb{E}\left[Y_n \, \big| \, \mathcal{F}_{n-1}\right]
=(\alpha + \beta)^{-(n-1)} Y_{n-1} = M_{n-1},
\label{F8}
\end{equation}
i.e. $M_n$ is an $\mathcal{F}_n$-martingale with
\begin{equation}
\mathbb{E}\left[M_n\right] = \mathbb{E}\left[M_0\right] = \mathbb{E}\left[Y_0\right] = 1
\qquad
\text{for all } n \geq 0
\label{F9}
\end{equation}
and, consequently, in view of \eqref{F7},
\begin{equation}
\mathbb{E}\left[Y_n\right] = \mathbb{E}\left[\left(\alpha + \beta\right)^n M_n\right]
= \left(\alpha + \beta\right)^n
\label{F10}
\end{equation}
(which also follows from \eqref{F5}).

Since $M_n \geq 0$ we have that \cite{D}
\begin{equation}
M_n \to M \ \text{a.s.},
\label{F11}
\end{equation}
where $M$ is a nonnegative random variable satisfying
\begin{equation}
0 \leq \mathbb{E}[M] \leq \mathbb{E}\left[M_n\right] = 1.
\label{F12}
\end{equation}

Equation \eqref{F1} implies that $\alpha + \beta \leq 1$. In the case
\begin{equation}
\alpha + \beta < 1
\label{F13}
\end{equation}
formulas \eqref{F7}, \eqref{F11}, and \eqref{F10} imply that
\begin{equation}
Y_n \to 0\ \text{a.s.}
\qquad
\text{and}
\qquad
\mathbb{E}\left[Y_n\right]  \to 0
\label{F14}
\end{equation}
(in other words $Y_n$ converges to $0$ almost surely and in the $L_1$-sense).
Also, in view of \eqref{F10},
\begin{equation}
\mathbb{E}\left[\sum_{n=0}^{\infty} Y_n\right] = \sum_{n=0}^{\infty} \mathbb{E}\left[Y_n\right]
=\sum_{n=0}^{\infty} \left(\alpha + \beta\right)^n
=\frac{1}{1-\alpha-\beta},
\label{F15}
\end{equation}
where the first equality is justified by Tonelli's Theorem. Therefore,
\begin{equation}
\sum_{n=0}^{\infty} Y_n < \infty \ \text{a.s.}
\label{F16}
\end{equation}
which suggests that, if $\alpha + \beta = \gamma < 1$, the forest will be, essentially, spared.

\medskip

\textbf{Remark 3.} The random variable $Y_n$ takes values in $\mathbb{N}$. Therefore, the fact that,
assuming $\alpha + \beta < 1$, we have $Y_n \to 0$ a.s. implies that there is a subset $\Omega_0$ of 
$\Omega$ with $\mathbb{P}(\Omega_0) = 1$ with the property that for every $\omega \in \Omega_0$ there 
is an $n_0 = n_0(\omega)$ such that $Y_n(\omega) = 0$ for every $n \geq n_0$. In particular, 
\eqref{F16} follows from \eqref{F14}. Also, in view of \eqref{F7}, it follows that $M_n(\omega) = 0$, 
too, for every $n \geq n_0$. Therefore, \eqref{F11} becomes $M_n \to 0$ a.s., i.e. $M = 0$ a.s. 
However, since $\mathbb{E}\left[M_n\right] = 1$ for all $n \geq 0$, the sequence $M_n$, $n \geq 0$, 
is not uniformly integrable \cite{D} (of course, in view of \eqref{F14}, $Y_n$, $n \geq 0$, 
is uniformly integrable). Finally, let us notice that $Y_n \to 0$ a.s. implies
\begin{equation}
\lim_n \mathbb{P}\left\{X_j^n = 0,\ 0 \leq j \leq n \right\}
= \lim_n \mathbb{P}\left\{Y_n = 0\right\} = 1.
\label{F17}
\end{equation}
\hfill $\diamondsuit$

\medskip

\textbf{Remark 4.} Suppose $\alpha + \beta < \gamma < 1$. Then, we can choose
$\tilde{\alpha} \geq \alpha$ and $\tilde{\beta} \geq \beta$ so that 
$\tilde{\alpha} + \tilde{\beta} = \gamma < 1$. Thus, by Remark 1 we get that formula \eqref{F14}
(and, consequently, \eqref{F16} thanks to the previous remark)
continues to hold, while \eqref{F15} becomes an inequality, namely
\begin{equation*}
\sum_{n=0}^{\infty} \mathbb{E}\left[Y_n\right]
\leq \frac{1}{1-\tilde{\alpha} - \tilde{\beta}}.
\end{equation*}
The same conclusions hold in the case where $\gamma < \alpha + \beta < 1$
(since we can choose $\tilde{\gamma} = \alpha + \beta$).
\hfill $\diamondsuit$

\subsection{The case $\alpha + \beta = \gamma = 1$}
In the case $\alpha + \beta = \gamma = 1$ formula \eqref{F7} becomes
\begin{equation}
M_n = Y_n
\label{E1}
\end{equation}
and, therefore $Y_n$ is an $\mathcal{F}_n$-martingale with
\begin{equation}
\mathbb{E}\left[Y_n\right] = \mathbb{E}\left[Y_0\right] = 1
\qquad
\text{for all } n \geq 0.
\label{E2}
\end{equation}
Formula \eqref{F15} now becomes
\begin{equation}
\mathbb{E}\left[\sum_{n=0}^{\infty} Y_n\right] = \infty
\label{E3}
\end{equation}
(hence, the fire might cause a serious damage to the forest).

Formula \eqref{F11} is, of course, still valid, namely
\begin{equation}
Y_n \to Y \text{a.s.},
\label{E4}
\end{equation}
where here $Y$ is a nonnegative random variable taking values in $\mathbb{N}$ (since $Y_n$
takes values in $\mathbb{N}$ for all $n \geq 0$) and satisfying
\begin{equation}
0 \leq \mathbb{E}[Y] \leq \mathbb{E}\big[Y_n\big] = 1.
\label{E5}
\end{equation}

\medskip

\textbf{Theorem 1.} If $\gamma = \alpha + \beta = 1$ (with $\alpha \beta > 0$, as
usual) and $Y_n$ is as in \eqref{C1}, then,
\begin{equation}
Y_n \to 0 \ \text{a.s.}
\label{E6}
\end{equation}

\smallskip

\textit{Proof}. From \eqref{E4} we know that $Y_n \to Y$ \text{a.s.}, i.e. that the almost sure
limit of $Y_n$ exists. Hence, we only need to show that $Y = 0$ \text{a.s.}

Consider the events
\begin{equation}
A_n = \left\{X_j^n = 1, \ 0 \leq j \leq n \right\},
\qquad
n \geq 0
\label{E7}
\end{equation}
(notice that $A_0 = \Omega$).

With the help of Example 1 and \eqref{F5} it is not hard to see that
\begin{equation}
\mathbb{P}\left(A_n\right) > 0
\qquad
\text{for all }
n \geq 0.
\label{E8}
\end{equation}

For notational convenience we denote by $\mathbb{P}_n$ the conditional probability given $A_n$,
that is for any event $B$ we have
\begin{equation}
\mathbb{P}_n(B)
= \mathbb{P}(B|A_n),
\qquad
n \geq 0
\label{E8a}
\end{equation}
(of course, $\mathbb{P}_0 = \mathbb{P}$).

Now, let
\begin{equation}
x_n = \mathbb{P}_n\left\{Y = 0\right\} = \mathbb{P}\left\{Y = 0 \,\big| \, A_n\right\}
= \mathbb{P}\left\{Y = 0 \,\big| \, X_j^n = 1, \ 0 \leq j \leq n\right\},
\quad
n \geq 0.
\label{E9}
\end{equation}

Then, we have
\begin{align}
x_0 &= \mathbb{P}\left\{Y = 0\right\}
\nonumber
\\
&= \sum_{\epsilon_0, \epsilon_1 = 0  \text{ or } 1}
\mathbb{P}\left\{Y=0 \, \big| \, X_0^1 = \epsilon_0, \  X_1^1 = \epsilon_1 \right\}
\mathbb{P}\left\{X_0^1 = \epsilon_0, \  X_1^1 = \epsilon_1 \right\}
\nonumber
\\
&= \sum_{\epsilon_0, \epsilon_1 = 0 \text{ or } 1}
\mathbb{P}\left\{Y=0 \, \big| \, X_0^1 = \epsilon_0, \  X_1^1 = \epsilon_1 \right\}
\mathbb{P}\left\{X_0^1 = \epsilon_0 \right\}
\mathbb{P}\left\{X_1^1 = \epsilon_1 \right\}
\label{E10}
\end{align}
(the second equality follows from \eqref{B2} and the fact that $\mathcal{G}_0$ is the trivial
$\sigma$-algebra). Notice that the sum in \eqref{E10} consists of four terms.

The next three formulas follow from \eqref{A1} and \eqref{E9}.
\begin{equation}
\mathbb{P}\left\{Y=0 \, \big| \, X_0^1 = 0, \  X_1^1 = 0 \right\} = 1,
\label{E11a}
\end{equation}
\begin{equation}
\mathbb{P}\left\{Y=0 \, \big| \, X_0^1 = 0, \  X_1^1 = 1 \right\}
= \mathbb{P}\left\{Y=0 \, \big| \, X_0^1 = 1, \  X_1^1 = 0 \right\} = x_0
\label{E11b}
\end{equation}
(both probabilities in \eqref{E11b} are equal to the conditional probability of the event $\{Y=0\}$ given 
that the fire started at one point, which is exactly the probability denoted by $x_0$),
\begin{equation}
\mathbb{P}\left\{Y=0 \, \big| \, X_0^1 = 1, \  X_1^1 = 1 \right\} = x_1.
\label{E11c}
\end{equation}
Substituting \eqref{E11a}, \eqref{E11b}, and \eqref{E11c} in \eqref{E10}
yields
\begin{align}
x_0 &= \mathbb{P}\left\{X_0^1 = 0 \right\} \mathbb{P}\left\{X_1^1 = 0 \right\}
+ x_0 \mathbb{P}\left\{X_0^1 = 0 \right\} \mathbb{P}\left\{X_1^1 = 1 \right\}
\nonumber
\\
&+ x_0 \mathbb{P}\left\{X_0^1 = 1 \right\} \mathbb{P}\left\{X_1^1 = 0 \right\}
+ x_1 \mathbb{P}\left\{X_0^1 = 1 \right\} \mathbb{P}\left\{X_1^1 = 1 \right\}.
\label{E12}
\end{align}
The above probabilities can be computed by \eqref{F5} (and the assumption that $\alpha + \beta = 1$).
We have
\begin{equation}
\mathbb{P}\left\{X_0^1 = 0 \right\} = 1 - \beta = \alpha,
\qquad
\mathbb{P}\left\{X_0^1 = 1 \right\} = \beta,
\label{E12a}
\end{equation}
and
\begin{equation}
\mathbb{P}\left\{X_1^1 = 0 \right\} = 1 - \alpha = \beta,
\qquad
\mathbb{P}\left\{X_0^1 = 1 \right\} = \alpha.
\label{E12b}
\end{equation}
Hence, by using \eqref{E12a} and \eqref{E12b} in \eqref{E12} we get
\begin{equation}
x_0 = \alpha \beta + \alpha^2 x_0 + \beta^2 x_0 + \alpha \beta x_1.
\label{E13}
\end{equation}
Using the fact that $\alpha^2 + \beta^2 = 1 - 2 \alpha \beta$, formula \eqref{E13} yields
\begin{equation}
x_0 = \frac{1 + x_1}{2} > \frac{1}{2}.
\label{E14}
\end{equation}

Next, we will look for a recursive formula for $x_n$.
\begin{align}
x_n &= \mathbb{P}_n\left\{Y = 0\right\}
\nonumber
\\
&= \sum_{\epsilon_j = 0 \text{ or } 1 \atop 0 \leq j \leq n+1}
\left[\,\mathbb{P}_n\left\{Y=0 \, \big| \, X_j^{n+1} = \epsilon_j, \ 0 \leq j \leq n+1 \right\} \right.
\nonumber
\\
&\qquad \qquad \ \ \left. \times\,
\mathbb{P}_n\left\{X_j^{n+1} = \epsilon_j, \ 0 \leq j \leq n+1 \right\}\right]
\label{E15}
\end{align}
From \eqref{E7} and \eqref{B1} we have that $A_n \in \mathcal{G}_n$. Hence, in view of \eqref{B2} and \eqref{E8a}, the events $\{X_j^{n+1} = \epsilon_j\}$, $0 \leq j \leq n+1$, are independent with respect to 
the probability $\mathbb{P}_n$ and, consequently, formula \eqref{E15} becomes
\begin{align}
x_n &= \mathbb{P}_n\left\{Y = 0\right\}
\nonumber
\\
&= \sum_{\epsilon_j = 0 \text{ or } 1 \atop 0 \leq j \leq n+1}
\left[\,\mathbb{P}_n\left\{Y=0 \, \big| \, X_j^{n+1} = \epsilon_j, \ 0 \leq j \leq n+1 \right\} \right.
\nonumber
\\
&\qquad \qquad \ \  \times\,
\prod_{j=0}^{n+1}\mathbb{P}_n\left\{X_j^{n+1} = \epsilon_j\right\}\big].
\label{E16}
\end{align}
Now, in view of \eqref{A1}, the assumption that $\gamma = 1$ implies
(recall \eqref{E7} and \eqref{E8a})
\begin{equation}
\mathbb{P}_n\left\{X_j^{n+1} = 0\right\} = 0
\,\, \Longleftrightarrow \,\,
\mathbb{P}_n\left\{X_j^{n+1} = 1\right\} = 1
\qquad
\text{for every } j \in \{1,2, \ldots, n\},
\label{E17}
\end{equation}
thus, by using \eqref{E17} and \eqref{E16} we obtain
\begin{align}
x_n &= \mathbb{P}_n\left\{Y = 0\right\}
\nonumber
\\
&= \sum_{\epsilon_0, \epsilon_1 = 0 \text{ or } 1}
\left[\,\mathbb{P}_n\left\{Y=0 \, \big| \, 
X_0^{n+1} = \epsilon_0, \ X_{n+1}^{n+1} = \epsilon_1, \ 
X_j^{n+1} = 1, \ 1 \leq j \leq n\right\} \right.
\nonumber
\\
&\qquad \qquad \qquad \ \ \left. \times\,
\mathbb{P}_n\left\{X_0^{n+1} = \epsilon_0\right\}
\mathbb{P}_n\left\{X_{n+1}^{n+1} = \epsilon_1\right\}\right]
\nonumber
\\
&= \sum_{\epsilon_0, \epsilon_1 = 0 \text{ or } 1}
\left[\,\mathbb{P}\left\{Y=0 \, \big| \, 
X_0^{n+1} = \epsilon_0, \ X_{n+1}^{n+1} = \epsilon_1, \ 
X_j^{n+1} = 1, \ 1 \leq j \leq n\right\} \right.
\nonumber
\\
&\qquad \qquad \qquad \ \ \left. \times\,
\mathbb{P}_n\left\{X_0^{n+1} = \epsilon_0\right\}
\mathbb{P}_n\left\{X_{n+1}^{n+1} = \epsilon_1\right\}\right],
\label{E18}
\end{align}
where the last equality (i.e, to use $\mathbb{P}$ instead of $\mathbb{P}_n$ in the conditional
probability) follows from the first equality of \eqref{B2} and the fact that (in view of \eqref{B1})
\begin{equation*}
\left\{X_0^{n+1} = \epsilon_0, \ X_{n+1}^{n+1} = \epsilon_1, \ X_j^{n+1} = 1, \ 1 \leq j \leq n\right\}
\in \mathcal{G}_{n+1}.
\end{equation*}

In view of \eqref{E9}, we have for $n \geq 1$
\begin{equation}
\mathbb{P}\left\{Y=0 \, \big| \, X_0^{n+1} = 0, \ X_{n+1}^{n+1} = 0, \ 
X_j^{n+1} = 1, \ 1 \leq j \leq n\right\} = x_{n-1},
\label{E19a}
\end{equation}
\begin{equation}
\mathbb{P}\left\{Y=0 \, \big| \, X_0^{n+1} = 0, \ X_{n+1}^{n+1} = 1, \ 
X_j^{n+1} = 1, \ 1 \leq j \leq n\right\} = x_n,
\label{E19b}
\end{equation}
\begin{equation}
\mathbb{P}\left\{Y=0 \, \big| \, X_0^{n+1} = 1, \ X_{n+1}^{n+1} = 0, \ 
X_j^{n+1} = 1, \ 1 \leq j \leq n\right\} = x_n,
\label{E19c}
\end{equation}
and
\begin{equation}
\mathbb{P}\left\{Y=0 \, \big| \, X_0^{n+1} = 1, \ X_{n+1}^{n+1} = 1, \ 
X_j^{n+1} = 1, \ 1 \leq j \leq n\right\} = x_{n+1}.
\label{E19d}
\end{equation}
Substituting \eqref{E19a}, \eqref{E19b}, \eqref{E19c}, and \eqref{E19d} in \eqref{E18}
yields
\begin{align}
x_n &= x_{n-1}\mathbb{P}_n\left\{X_0^{n+1} = 0 \right\} \mathbb{P}_n\left\{X_{n+1}^{n+1} = 0 \right\}
+ x_n \mathbb{P}_n\left\{X_0^{n+1} = 0 \right\} \mathbb{P}_n\left\{X_{n+1}^{n+1} = 1 \right\}
\nonumber
\\
&+ x_n \mathbb{P}_n\left\{X_0^{n+1} = 1 \right\} \mathbb{P}_n\left\{X_{n+1}^{n+1} = 0 \right\}
+ x_{n+1} \mathbb{P}_n\left\{X_0^{n+1} = 1 \right\} \mathbb{P}_n\left\{X_{n+1}^{n+1} = 1 \right\}.
\label{E20}
\end{align}
Finally, we calculate the probabilities that appear in \eqref{E20}. In view of \eqref{E7}, \eqref{E8a},\eqref{B2}, \eqref{A1}, and \eqref{A4}, we get
\begin{align}
\mathbb{P}_n\left\{X_0^{n+1} = 0 \right\}
&= \mathbb{P}\left\{X_0^{n+1} = 0 \, \big| \, X_j^n = 1, \ 0 \leq j \leq n\right\}
\nonumber
\\
&= \mathbb{P}\left\{X_0^{n+1} = 0 \, \big| \, X_0^n = 1\right\}
= 1 - \beta = \alpha
\label{E21a}
\end{align}
and
\begin{align}
\mathbb{P}_n\left\{X_{n+1}^{n+1} = 0 \right\}
&= \mathbb{P}\left\{X_{n+1}^{n+1} = 0 \, \big| \, X_j^n = 1, \ 0 \leq j \leq n\right\}
\nonumber
\\
&= \mathbb{P}\left\{X_{n+1}^{n+1} = 0 \, \big| \, X_n^n = 1\right\}
= 1 - \alpha = \beta,
\label{E21b}
\end{align}
while from \eqref{E21a} and \eqref{E21b} we get immediately that
\begin{equation}
\mathbb{P}_n\left\{X_0^{n+1} = 1 \right\} = \beta
\qquad
\text{and}
\qquad
\mathbb{P}_n\left\{X_{n+1}^{n+1} = 1 \right\} = \alpha.
\label{E21c}
\end{equation}
Therefore, by substituting \eqref{E21a}, \eqref{E21b}, and \eqref{E21c}, in \eqref{E20} we obtain
\begin{equation}
x_n = \alpha \beta x_{n-1} + \alpha^2 x_n + \beta^2 x_n + \alpha \beta x_{n+1},
\label{E22}
\end{equation}
which implies (since $\alpha^2 + \beta^2 = 1 - 2\alpha \beta$)
\begin{equation}
x_n = \frac{x_{n-1} + x_{n+1}}{2}
\qquad
n \geq 1.
\label{E23}
\end{equation}
Equation \eqref{E23} is a very simple second-order linear difference equation with constant coefficients.
Its general solution is
\begin{equation*}
x_n = c_1 + c_2 n,
\qquad
n \geq 0,
\end{equation*}
where $c_1, c_2$ are constants.
However, in our case $x_n$ is a probability. Hence, we must have $c_2 = 0$ and, consequently,
$x_n$ is constant, i.e.
\begin{equation}
x_n = x_0
\qquad
n \geq 1.
\label{E24}
\end{equation}
In particular, $x_1 = x_0$, hence \eqref{E14} becomes
\begin{equation*}
x_0 = \frac{1 + x_0}{2} 
\end{equation*}
which implies that $x_0 = 1$, i.e. (in view of \eqref{E9})
$\mathbb{P}\{Y = 0\} = 1$.
\hfill $\blacksquare$

\medskip

The symbol $\blacksquare$ indicates the end of a proof.

\medskip

Actually, in the proof of Theorem 1 we have shown that
\begin{equation}
x_n = \mathbb{P}\left\{Y=0 \, \big| \, X_j^n = 1, \ 0 \leq j \leq n\right\} = 1
\qquad
\text{for all }
n \geq 0,
\label{E25}
\end{equation}
which is stronger than $\mathbb{P}\{Y = 0\} = 1$.

\medskip

\textbf{Remark 5.} From Theorem 1 we see that $\mathbb{E}[Y] = 0 \ne 1 = \mathbb{E}[Y_n]$. Hence,
in the case $\gamma = \alpha + \beta = 1$ the sequence $Y_n$, $n \geq 0$,
is not uniformly integrable \cite{D}.
\hfill $\diamondsuit$

\medskip

\textbf{Remark 6.} Suppose $\alpha + \beta = 1$ and $\gamma < 1$. Then, from formula \eqref{C2}
we get that
\begin{equation*}
\mathbb{E}\left[Y_n \, \big| \, \mathcal{F}_{n-1}\right]
\leq Y_{n-1},
\end{equation*}
i.e. $Y_n$ is a (nonnegative) $\mathcal{F}_n$-supermartingale. Therefore, \eqref{E4}
continues to hold \cite{D}. Consequently, in view of Remark 1, we get that Theorem 1 and, furthermore, 
\eqref{E25} are still valid.
\hfill $\diamondsuit$

\medskip

Regarding the fate of the forest, the results of this subsection lead us to the conclusion that
$\alpha + \beta = \gamma = 1$ is a borderline case. On one hand formula \eqref{E3} suggests that
many trees will be burnt, while on the other hand formula \eqref{E25} tells us that eventually the fire 
will (eventually) die out.

\subsection{The growth of $\mathbb{E}\left[Y_n^2\right]$}
In this subsection we focus on $Y_n^2$. Actually, our main interest is the behavior of
$\mathbb{E}\left[Y_n^2\right]$ as $n \to \infty$.

From the defining formula \eqref{C1} of $Y_n$ we get that
\begin{align}
Y_n^2 = \left(\sum_{j=0}^n X_j^n\right)^2 &= \sum_{j=0}^n \left(X_j^n\right)^2
+ 2\sum_{j=0}^n \sum_{k=j+1}^n X_j^n X_k^n
\nonumber
\\
&= \sum_{j=0}^n X_j^n
+ 2\sum_{j=0}^n \sum_{k=j+1}^n X_j^n X_k^n
\nonumber
\\
&= Y_n + 2\sum_{j=0}^n \sum_{k=j+1}^n X_j^n X_k^n,
\qquad\qquad
n \geq 0.
\label{G1}
\end{align}
Hence,
\begin{equation}
\mathbb{E}\left[Y_n^2\right] = 
\mathbb{E}\left[Y_n\right] + 2\sum_{j=0}^n \sum_{k=j+1}^n \mathbb{E}\left[X_j^n X_k^n\right],
\qquad
n \geq 0.
\label{G2}
\end{equation}
Throughout this section we assume \eqref{F1}, i.e. $\gamma = \alpha + \beta$. Therefore, \eqref{F10} 
holds and, consequently, \eqref{G2} becomes
\begin{equation}
\mathbb{E}\left[Y_n^2\right] = 
\gamma^n + 2\sum_{j=0}^n \sum_{k=j+1}^n \mathbb{E}\left[X_j^n X_k^n\right],
\qquad
n \geq 0.
\label{G3}
\end{equation}

Let us set
\begin{equation}
R_0 = 0,
\qquad
R_n = \sum_{j=0}^n \sum_{k=j+1}^n \mathbb{E}\left[X_j^n X_k^n\right],
\qquad
n \geq 1,
\label{G4}
\end{equation}
so that \eqref{G3} can be written as
\begin{equation}
\mathbb{E}\left[Y_n^2\right] = 
\gamma^n + 2 R_n,
\qquad
n \geq 0.
\label{G5}
\end{equation}
Now, under \eqref{F1} formula \eqref{B2e} becomes
\begin{align}
&\mathbb{E}\left[X_j^n X_k^n \, \big| \, \mathcal{F}_{n-1}\right]
=\left(\alpha X_{j-1}^{n-1} + \beta X_j^{n-1} \right)
\left(\alpha X_{k-1}^{n-1} + \beta X_k^{n-1} \right)
\nonumber
\\
&=\alpha^2 X_{j-1}^{n-1} X_{k-1}^{n-1} + \alpha \beta X_{j-1}^{n-1} X_k^{n-1}
+ \alpha \beta X_j^{n-1} X_{k-1}^{n-1} + \beta^2 X_j^{n-1} X_k^{n-1},
\label{G6}
\end{align}
hence, by taking expectations we obtain
\begin{align}
\mathbb{E}\left[X_j^n X_k^n \right] &= \alpha^2 \mathbb{E}\left[X_{j-1}^{n-1} X_{k-1}^{n-1} \right]
+ \alpha \beta \mathbb{E}\left[X_{j-1}^{n-1} X_k^{n-1} \right]
\nonumber
\\
&+ \alpha \beta \mathbb{E}\left[X_j^{n-1} X_{k-1}^{n-1} \right]
+ \beta^2 \mathbb{E}\left[X_j^{n-1} X_k^{n-1} \right].
\label{G7}
\end{align}
Using \eqref{G7} in \eqref{G4} yields
\begin{align}
R_n &= \alpha^2 \sum_{j=0}^n \sum_{k=j+1}^n \mathbb{E}\left[X_{j-1}^{n-1} X_{k-1}^{n-1} \right]
+ \alpha \beta \sum_{j=0}^n \sum_{k=j+1}^n \mathbb{E}\left[X_{j-1}^{n-1} X_k^{n-1} \right]
\nonumber
\\
&+ \alpha \beta \sum_{j=0}^n \sum_{k=j+1}^n \mathbb{E}\left[X_j^{n-1} X_{k-1}^{n-1} \right]
+ \beta^2 \sum_{j=0}^n \sum_{k=j+1}^n \mathbb{E}\left[X_j^{n-1} X_k^{n-1} \right].
\label{G8}
\end{align}
In view of the boundary conditions \eqref{A4}, formula \eqref{G8} becomes
\begin{align}
R_n &= 
\alpha^2 \sum_{j=0}^{n-1} \sum_{k=j+1}^{n-1} \mathbb{E}\left[X_j^{n-1} X_k^{n-1} \right]
+ \alpha \beta \sum_{j=0}^{n-1} \sum_{k=j+2}^{n-1} \mathbb{E}\left[X_j^{n-1} X_k^{n-1} \right]
\nonumber
\\
&+ \alpha \beta \sum_{j=0}^{n-1} \sum_{k=j}^{n-1} \mathbb{E}\left[X_j^{n-1} X_k^{n-1} \right]
+ \beta^2 \sum_{j=0}^{n-1} \sum_{k=j+1}^{n-1} \mathbb{E}\left[X_j^{n-1} X_k^{n-1} \right]
\nonumber
\\
&= (\alpha^2 + \beta^2) R_{n-1}
\nonumber
\\
&+ \alpha \beta \sum_{j=0}^{n-1} \sum_{k=j+2}^{n-1} \mathbb{E}\left[X_j^{n-1} X_k^{n-1} \right]
+ \alpha \beta \sum_{j=0}^{n-1} \sum_{k=j}^{n-1} \mathbb{E}\left[X_j^{n-1} X_k^{n-1} \right]
\nonumber
\\
&=(\alpha^2 + \beta^2) R_{n-1}
\nonumber
\\
&+ \alpha \beta \left(R_{n-1} - \sum_{j=0}^{n-1}  \mathbb{E}\left[X_j^{n-1} X_{j+1}^{n-1} \right]\right)
+ \alpha \beta \left(R_{n-1} + \sum_{j=0}^{n-1} \mathbb{E}\left[X_j^{n-1} X_j^{n-1} \right]\right)
\nonumber
\\
&=(\alpha + \beta)^2 R_{n-1}
- \alpha \beta \sum_{j=0}^{n-2}  \mathbb{E}\left[X_j^{n-1} X_{j+1}^{n-1} \right]
+ \alpha \beta\sum_{j=0}^{n-1} \mathbb{E}\left[X_j^{n-1} \right]
\nonumber
\\
&=(\alpha + \beta)^2 R_{n-1}
- \alpha \beta \sum_{j=0}^{n-2}  \mathbb{E}\left[X_j^{n-1} X_{j+1}^{n-1} \right]
+ \alpha \beta \mathbb{E}\left[Y_{n-1} \right],
\label{G9}
\end{align}
hence, in view of \eqref{F1} and \eqref{F10},
\begin{equation}
R_n = \gamma^2 R_{n-1} + \alpha \beta \gamma^{n-1}
- \alpha \beta \sum_{j=0}^{n-2}  \mathbb{E}\left[X_j^{n-1} X_{j+1}^{n-1} \right],
\qquad
n \geq 1.
\label{G10}
\end{equation}
One very simple implication of formula \eqref{G10} is the inequality
\begin{equation}
R_n \leq \gamma^2 R_{n-1} + \alpha \beta \gamma^{n-1},
\qquad
n \geq 1.
\label{G11}
\end{equation}

\medskip

\underline{\textbf{Case} $\gamma = \alpha + \beta < 1$}\textbf{.} In this case \eqref{G11} implies
\begin{equation*}
\gamma^{-2n} R_n - \gamma^{-2(n-1)} R_{n-1}  \leq \alpha \beta \gamma^{-(n+1)}
\qquad
n \geq 1
\end{equation*}
hence,
\begin{equation*}
\sum_{n=1}^N \left(\gamma^{-2n} R_n - \gamma^{-2(n-1)} R_{n-1}\right)
\leq \alpha \beta \sum_{n=1}^N \gamma^{-(n+1)},
\end{equation*}
i.e.
\begin{equation*}
\gamma^{-2N} R_N
\leq \alpha \beta \sum_{n=1}^N \gamma^{-(n+1)}
= \frac{\alpha \beta}{\gamma} \frac{\gamma^{-N} - 1}{1 - \gamma},
\qquad
N \geq 1,
\end{equation*}
from which it follows that
\begin{equation}
R_n \leq \alpha \beta \frac{1 - \gamma^n}{1 - \gamma} \gamma^{n-1},
\qquad
n \geq 1.
\label{G12}
\end{equation}

Finally, using \eqref{G12} in \eqref{G5} yields
\begin{equation}
\mathbb{E}\left[Y_n^2\right] \leq 
\left(\gamma + 2 \alpha \beta \frac{1 - \gamma^n}{1 - \gamma}\right) \gamma^{n-1},
\qquad
n \geq 1,
\label{G13}
\end{equation}
hence if $\gamma = \alpha + \beta < 1$ (more generally, if $\gamma \leq \alpha + \beta < 1$), then, 
as $n \to \infty$, $\mathbb{E}\left[Y_n^2\right]$ aproaches $0$ exponentially fast. In particular,
$Y_n \to 0$ in the $L_2$-sense.

\medskip

\underline{\textbf{Case} $\gamma = \alpha + \beta = 1$}\textbf{.} Here the fact that
(recall \eqref{F8} and \eqref{E1}) $Y_n$ is an $\mathcal{F}_n$-martingale
implies \cite{D} that $Y_n^2$ is an $\mathcal{F}_n$-submartingale. In particular,
\begin{equation}
\mathbb{E}\left[Y_n^2\right] \geq \mathbb{E}\left[Y_{n-1}^2\right],
\qquad
n \geq 1.
\label{G14}
\end{equation}
Then, the fact that the sequence $Y_n$, $n \geq 0$,
is not uniformly integrable (Remark 5) implies that
\begin{equation}
\lim_n \mathbb{E}\left[Y_n^2\right] = \infty,
\label{G14a}
\end{equation}
since, if $\mathbb{E}\left[Y_n^2\right]$ were bounded, then $Y_n$, $n \geq 0$,
would have been uniformly integrable \cite{D}.

We can easily get an upper estimate of the growth of $\mathbb{E}\left[Y_n^2\right]$.
From \eqref{G5} (since here $\gamma = 1$) we get that
\begin{equation}
\mathbb{E}\left[Y_n^2\right] = 1 + 2 R_n,
\qquad
n \geq 0.
\label{G15}
\end{equation}
Now, formula \eqref{G10} becomes
\begin{equation}
R_n = R_{n-1} + \alpha \beta 
- \alpha \beta \sum_{j=0}^{n-2}  \mathbb{E}\left[X_j^{n-1} X_{j+1}^{n-1} \right],
\qquad
n \geq 1,
\label{G16}
\end{equation}
while inequality \eqref{G11} becomes
\begin{equation}
R_n \leq R_{n-1} + \alpha \beta,
\qquad
n \geq 1.
\label{G17}
\end{equation}
Thus, in view of \eqref{G15}, the inequality \eqref{G17} implies immediately the upper 
estimate
\begin{equation}
\mathbb{E}\left[Y_n^2\right] \leq 1 + 2 \alpha \beta n,
\qquad
n \geq 0,
\label{G18}
\end{equation}
hence $\mathbb{E}\left[Y_n^2\right]$ grows at most linearly.

\section{The case $\alpha + \beta > 1$}
In the case where $\alpha + \beta > 1$ we have not established the existence of the almost sure limit
of $Y_n$. However, if we set
\begin{equation}
W_n = {\bf 1}_{\{Y_n = 0\}},
\qquad
n \geq 0,
\label{H1}
\end{equation}
then formula \eqref{C1b} implies
\begin{equation}
W_n \leq W_{n+1} \text{ a.s.},
\qquad
n \geq 0,
\label{H2}
\end{equation}
hence,
\begin{equation}
W_n \to W \text{ a.s.},
\label{H3}
\end{equation}
where $W$ is a Bernoulli random variable. Thus, by the Monotone Convergence Theorem we obtain that
(since $\mathbb{P}\left\{Y_n = 0 \right\} = \mathbb{E}\left[W_n\right]$)
\begin{equation}
\lim_n \mathbb{P}\left\{Y_n = 0 \right\} = \lim_n \mathbb{E}\left[W_n\right]
= \mathbb{E}[W].
\label{H4}
\end{equation}

If $\alpha = 1$, then the definition of the random field $X_j^n$ implies that $X_n^n = 1$ a.s. and,
consequently, $\mathbb{P}\left\{Y_n = 0 \right\} = 0$ for all $n \geq 0$, hence $W = 0$ a.s.
Likewise, if $\beta = 1$, then $X_0^n = 1$ a.s. hence, again
$\mathbb{P}\left\{Y_n = 0 \right\} = 0$ for all $n \geq 0$ and $W = 0$ a.s. For this reason, from now on
we will assume that
\begin{equation}
\alpha, \beta \ne 1.
\label{H9}
\end{equation}
Then, from the definition of the random field $X_j^n$ it is not hard to see that
\begin{equation}
\mathbb{P}\left\{Y_n = 0 \right\} < \mathbb{P}\left\{Y_{n+1} = 0 \right\},
\qquad
n \geq 0.
\label{H5}
\end{equation}

\subsection{The case $\alpha + \beta > 1$ and $\gamma = 1$}
Before giving the first result of this subsection let us notice that the assumption that
$\alpha + \beta > 1$ implies that $(1 - \alpha)(1 - \beta) < \alpha \beta$ and, hence,
\begin{equation}
\frac{(1 - \alpha)(1 - \beta)}{\alpha \beta} < 1.
\label{H12α}
\end{equation}

The first result of the subsection is a lower bound of the limit $\lim_n\mathbb{P}\left\{Y_n = 0 \right\}$.

\medskip

\textbf{Theorem 2.} If $\alpha + \beta > 1$ and $\gamma = 1$ (with $\alpha, \beta \ne 1$), while $Y_n$
is as in \eqref{C1}, then
\begin{equation}
\lim_n\mathbb{P}\left\{Y_n = 0 \right\} \geq \frac{(1 - \alpha)(1 - \beta)}{\alpha \beta}.
\label{H17}
\end{equation}

\smallskip

\textit{Proof}. Imitating the proof of Theorem 1 we set (recall \eqref{E7}, \eqref{E8a}, and \eqref{E9})
\begin{equation}
x_n = \mathbb{E}_n[W] = \mathbb{E}\left[W \,\big| \, A_n\right]
= \mathbb{E}\left[W \,\big| \, X_j^n = 1, \ 0 \leq j \leq n\right],
\qquad
n \geq 0
\label{H6}
\end{equation}
($\mathbb{E}_n$ is the expectation corresponding to the probability $\mathbb{P}_n$).

It is clear from \eqref{H6} that
\begin{equation}
x_n \geq x_{n+1},
\qquad
n \geq 0.
\label{H6b}
\end{equation}
In view of \eqref{H4} (and the fact that $\mathbb{E}_0 = \mathbb{E}$),
\begin{equation}
x_0 = \mathbb{E}[W] = \lim_n\mathbb{P}\left\{Y_n = 0 \right\}.
\label{H7}
\end{equation}
Also, in view of \eqref{H5} and \eqref{H6} we have that
\begin{align}
x_0 > \mathbb{P}\left\{Y_1 = 0 \right\} &= \mathbb{P}\left\{X_0^1 = X_1^1 = 0 \right\}
\nonumber
\\
&= \mathbb{P}\left\{X_0^1 = 0 \right\} \mathbb{P}\left\{X_1^1 = 0 \right\} = (1 - \beta)(1 - \alpha)
\label{H8}
\end{align}
(the second equality follows from \eqref{B2} and the fact that $\mathcal{G}_0$ is the trivial
$\sigma$-algebra; the third equality follows from \eqref{F5}).

At this point we notice that (with $\gamma = 1$) all formulas in the proof of Theorem 1 remain valid as 
long as they do not depend on the relation $\alpha + \beta = 1$. In particular, formula \eqref{E12} yields
\begin{equation*}
x_0 = (1 - \beta)(1 - \alpha) + [(1 - \beta)\alpha + \beta (1 - \alpha)] x_0 + \alpha \beta x_1,
\end{equation*}
or, equivalently,
\begin{equation}
[(1 - \alpha)(1 - \beta) + \alpha \beta] x_0
= (1 - \alpha)(1 - \beta) + \alpha \beta x_1,
\label{H10}
\end{equation}
while formula \eqref{E20} yields
\begin{equation*}
x_n
= (1 - \beta)(1 - \alpha) x_{n-1} + [(1 - \beta)\alpha + \beta (1 - \alpha)] x_n + \alpha \beta x_{n+1},
\qquad
n \geq 1,
\end{equation*}
or, equivalently,
\begin{equation}
[(1 - \alpha)(1 - \beta) + \alpha \beta] x_n
= (1 - \alpha)(1 - \beta) x_{n-1} + \alpha \beta x_{n+1},
\qquad
n \geq 1.
\label{H11}
\end{equation}
Equation \eqref{H11} is a very simple second-order linear difference equation with constant coefficients.
Its general solution is
\begin{equation}
x_n = c_1 + c_2 \left[\frac{(1 - \alpha)(1 - \beta)}{\alpha \beta}\right]^n,
\qquad
n \geq 0,
\label{H12}
\end{equation}
where $c_1, c_2$ are constants.

Using \eqref{H12} in \eqref{H10} yields
\begin{equation*}
[(1 - \alpha)(1 - \beta) + \alpha \beta] (c_1 + c_2)
= (1 - \alpha)(1 - \beta)
+ \alpha \beta \left[c_1 + c_2 \frac{(1 - \alpha)(1 - \beta)}{\alpha \beta} \right],
\end{equation*}
which implies that
\begin{equation}
c_2 = \frac{(1 - \alpha)(1 - \beta)}{\alpha \beta} (1 -c_1).
\label{H13}
\end{equation}
Thus, \eqref{H12} becomes
\begin{equation}
x_n = c_1 + (1 - c_1) \left[\frac{(1 - \alpha)(1 - \beta)}{\alpha \beta}\right]^{n+1},
\qquad
n \geq 0,
\label{H14}
\end{equation}
and the fact that $x_n$ is a decreasing sequence of probabilities implies (in view of \eqref{H12α})
that
\begin{equation}
0 \leq c_1 \leq 1.
\label{H15}
\end{equation}
Finally, for $n=0$ formula \eqref{H14} yields
\begin{equation}
\lim_n\mathbb{P}\left\{Y_n = 0 \right\} = x_0 = c_1 + (1 - c_1) \frac{(1 - \alpha)(1 - \beta)}{\alpha \beta},
\label{H16}
\end{equation}
from which \eqref{H17} follows.
\hfill $\blacksquare$

\medskip

\textbf{Remark 7.} The inequality \eqref{H17} gives a lower bound for 
$\lim_n\mathbb{P}\left\{Y_n = 0 \right\}$. In the extreme case where $\alpha = 1$ or $\beta = 1$ we have
seen that $\mathbb{P}\left\{Y_n = 0 \right\} = 0$ for all $n \geq 0$, hence \eqref{H17} becomes equality.
Also, in the case $\alpha + \beta = 1$ the right-hand side of \eqref{H17} is equal to $1$, hence, again
\eqref{H17} becomes equality. It is tempting to conjecture that \eqref{H17} becomes equality whenever
$\alpha + \beta > 1$ and $\gamma = 1$.
\hfill $\diamondsuit$

\medskip

The next theorem gives upper and lower bounds for the expectation $\mathbb{E}\left[Y_n \right]$.

\medskip

\textbf{Theorem 3.} If $\alpha + \beta > 1$ and $\gamma = 1$ (with $\alpha, \beta \ne 1$), while $Y_n$
is as in \eqref{C1}, then:

(i) We have the upper bound
\begin{equation}
\mathbb{E}\left[Y_n \right] \geq 1 + (\alpha + \beta - 1)n,
\qquad
n \geq 0.
\label{H18}
\end{equation}

(ii) For any $\varepsilon > 0$ there is a $n_0 = n_0(\varepsilon)$ such that
\begin{equation}
\mathbb{E}\left[Y_n \right] \leq \left(\frac{\alpha + \beta - 1}{\alpha \beta} + \varepsilon\right) n,
\qquad
n \geq n_0
\label{H23}
\end{equation}
(notice that $(\alpha + \beta - 1) / (\alpha \beta) < 1$), hence
\begin{equation}
\limsup_n \frac{\mathbb{E}\left[Y_n \right]}{n} \leq 
\frac{\alpha + \beta - 1}{\alpha \beta}.
\label{H23a}
\end{equation}

\smallskip

\textit{Proof}. (i) For $\gamma = 1$ formula \eqref{C2} becomes
\begin{equation}
\mathbb{E}\left[Y_n \, \big| \, \mathcal{F}_{n-1}\right]
= Y_{n-1}
+ (\alpha + \beta - 1) \left(Y_{n-1} - \sum_{j=1}^{n-1} X_{j-1}^{n-1} X_j^{n-1}\right),
\quad
n \geq 1.
\label{H19}
\end{equation}
Now, the fact that
\begin{equation*}
Y_{n-1} = \sum_{j=0}^{n-1} X_j^{n-1} = \sum_{j=1}^{n} X_{j-1}^{n-1}
\end{equation*}
implies (since the $X_j^{n-1}$'s are Bernoulli random variables)
\begin{equation}
Y_{n-1} \geq 1 + \sum_{j=1}^{n-1} X_{j-1}^{n-1} X_j^{n-1}.
\label{H20}
\end{equation}
Thus, \eqref{H19} implies
\begin{equation}
\mathbb{E}\left[Y_n \, \big| \, \mathcal{F}_{n-1}\right]
\geq Y_{n-1}
+ \alpha + \beta - 1,
\qquad
n \geq 1,
\label{H21}
\end{equation}
and, therefore, by taking expectations we get
\begin{equation}
\mathbb{E}\left[Y_n\right] \geq \mathbb{E}\left[Y_{n-1}\right] + \alpha + \beta - 1,
\qquad
n \geq 1
\label{H22}
\end{equation}
(also, since $\alpha + \beta > 1$, $Y_n$ is an $\mathcal{F}_n$-submartingale), and \eqref{H18}
follows immediately from \eqref{H22} and the fact that $Y_0 = 1$.

(ii) From \eqref{H17} we get that, given $\varepsilon_1 > 0$ there is a $n_1 = n_1(\varepsilon_1)$ such 
that
\begin{equation}
\mathbb{P}\left\{Y_n = 0 \right\} \geq \frac{(1 - \alpha)(1 - \beta)}{\alpha \beta} - \varepsilon_1,
\qquad
n \geq n_1.
\label{H24}
\end{equation}

Now, since $Y_n \in \{0,1, \ldots, n+1\}$,
\begin{equation}
\mathbb{E}\left[Y_n\right] \leq (n+1) \mathbb{P}\left\{Y_n \geq 1 \right\}
= (n+1) \left(1 -\mathbb{P}\left\{Y_n = 0 \right\} \right),
\label{H25}
\end{equation}
hence, in view of \eqref{H24},
\begin{align}
\mathbb{E}\left[Y_n\right] 
&\leq (n+1) \left(1 - \frac{(1 - \alpha)(1 - \beta)}{\alpha \beta} + \varepsilon_1 \right)
\nonumber
\\
&= \left(\frac{\alpha + \beta - 1}{\alpha \beta} + \varepsilon_1\right) (n+1),
\qquad
n \geq n_1,
\nonumber
\end{align}
from which we easily obtain \eqref{H23}.
\hfill $\blacksquare$

\medskip

From the inequality \eqref{H18} we can easily get an upper bound for the probability
$\mathbb{P}\left\{Y_n = 0 \right\}$ and its limit
$\lim_n\mathbb{P}\left\{Y_n = 0 \right\}$.

\medskip

\textbf{Corollary 1.} If $\alpha + \beta > 1$ and $\gamma = 1$ (with $\alpha, \beta \ne 1$), while $Y_n$
is as in \eqref{C1}, then
\begin{equation}
\mathbb{P}\left\{Y_n = 0 \right\} \leq \frac{n}{n+1} (2 - \alpha - \beta),
\qquad
n \geq 0,
\label{H27}
\end{equation}
hence
\begin{equation}
\lim_n\mathbb{P}\left\{Y_n = 0 \right\} \leq 2 - \alpha - \beta.
\label{H28}
\end{equation}

\smallskip

\textit{Proof}. Using \eqref{H18} in \eqref{H28} yields
\begin{equation}
1 + (\alpha + \beta - 1)n \leq (n+1) \left(1 -\mathbb{P}\left\{Y_n = 0 \right\} \right),
\label{H29}
\end{equation}
from which \eqref{H27} follows.
\hfill $\blacksquare$

\medskip

Formulas \eqref{H18} and \eqref{H28} suggest that, in the case where $\alpha + \beta > 1$ and $\gamma = 1$
the forest will be considerably damaged by the fire.

\medskip

\textbf{Remark 8.} By Remark 1 formulas \eqref{H17} and \eqref{H23} remain valid in the case $\gamma < 1$.
\hfill $\diamondsuit$

\subsection{The random variables $T_n$, $n \geq 1$}
We set
\begin{equation}
T_n = \frac{Y_n}{n},
\qquad
n \geq 1.
\label{J1}
\end{equation}

Recall that $Y_n \in \{0,1, \ldots, n+1\}$. Hence,
\begin{equation}
0 \leq T_n \leq  \frac{n+1}{n} \leq 2,
\qquad
n \geq 1
\label{J2}
\end{equation}
(thus, for example, $\mathbb{E}\left[T_n^2\right] \leq 4$ for all $n \geq 1$).

Formula \eqref{J2} implies \cite{D} that the sequence $\{T_n\}_{n \geq 1}$ is uniformly integrable and, 
also, that the sequence of the distribution functions of $T_n$, $n \geq 1$, is tight.

Let us observe that from formulas \eqref{H18} and \eqref{H23a} it follows immediately that
\begin{equation}
\alpha + \beta - 1 \leq \liminf_n \mathbb{E}\left[T_n\right] 
\leq \limsup_n \mathbb{E}\left[T_n\right]
\leq \frac{\alpha + \beta - 1}{\alpha \beta}.
\label{J4}
\end{equation}

\medskip

\textbf{Open Question.} Does $T_n$ converge, at least in distribution?

\medskip

Notice that if $T_n$ converges in distribution to a random variable $T$, then we must also have that
$\lim_n \mathbb{E}\left[T_n\right] = \mathbb{E}[T]$.
This follows from the fact that the distributional convergence of $T_n$ implies \cite{D} that
there are random variables $\tilde{T}_n$, $n \geq 1$, such that $\tilde{T}_n$ and $T_n$ have the same distribution (for every $n \geq 1$) and $\tilde{T}_n \to \tilde{T}$ a.s., where $\tilde{T}$ and $T$
have the same distribution. Furthermore, \eqref{J2} implies that $\{\tilde{T}_n\}_{n \geq 1}$, too, is
uniformly integrable, hence $\lim_n \mathbb{E}\left[T_n\right]
= \lim_n \mathbb{E}[\tilde{T}_n] = \mathbb{E}[\tilde{T}] = \mathbb{E}[T]$.

\medskip

Finally, let us introduce some random variables which may be helpful in the study of the random field
$X^n_j$. We set
\begin{equation}
L_0 = \inf\left\{n \geq 0 : X^n_0 = 0\right\}
\label{J5}
\end{equation}
and
\begin{equation}
L_m = \inf\left\{n \geq 0 : X^{n + L_0 + L_1 + \cdots + L_{m-1}}_m = 0\right\},
\qquad
m \geq 1;
\label{J6}
\end{equation}
also
\begin{equation}
K_0 = \inf\left\{n \geq 0 : X^n_n = 0\right\}
\label{J5b}
\end{equation}
and
\begin{equation}
K_m = \inf\left\{n \geq 0 : X^{n + K_0 + K_1 + \cdots + K_{m-1}}_{n-m} = 0\right\},
\qquad
m \geq 1.
\label{J6b}
\end{equation}
Observe that
\begin{equation}
X^{n + L_0 + L_1 + \cdots + L_m}_m = X^{n + K_0 + K_1 + \cdots + K_m}_{n-m} = 0
\qquad
\text{for all }
n \geq 0.
\label{J7}
\end{equation}

Each of the above random variables resembles the variable $L$ of the appendix. The boundary condition
$X^0_0 = 1$ implies that $\mathbb{P}\{L_0 = 0\} = \mathbb{P}\{K_0 = 0\} = 0$ and, hence, formula 
\eqref{D10} of the appendix
together with \eqref{A1} imply that
\begin{equation}
L_0 \text{ is a geometric random variable with parameter }1-\beta
\label{J8}
\end{equation}
and
\begin{equation}
K_0 \text{ is a geometric random variable with parameter }1-\alpha.
\label{J8b}
\end{equation}

Regarding $L_m$ and $K_m$, $m \geq 1$, formulas \eqref{A1}, \eqref{D10a}, and \eqref{D10} imply that
\begin{equation}
\mathbb{P}\left\{L_m = 0\right\} = \mathbb{P}\left\{X^{L_0 + L_1 + \cdots + L_{m-1}}_m = 0\right\}
\label{J9}
\end{equation}
and
\begin{equation}
\mathbb{P}\left\{L_m = k\right\}
= \mathbb{P}\left\{X^{L_0 + L_1 + \cdots + L_{m-1}}_m = 1\right\} (1-\beta) \beta^{k-1},
\qquad
k \geq 1,
\label{J10}
\end{equation}
while (symmetrically)
\begin{equation}
\mathbb{P}\left\{K_m = 0\right\}
 = \mathbb{P}\left\{X^{K_0 + K_1 + \cdots + K_{m-1}}_{n-m} = 0\right\}
\label{J9a}
\end{equation}
and
\begin{equation}
\mathbb{P}\left\{K_m = k\right\}
= \mathbb{P}\left\{X^{K_0 + K_1 + \cdots + K_{m-1}}_{n-m} = 1\right\} (1-\alpha) \alpha^{k-1},
\qquad
k \geq 1.
\label{J10a}
\end{equation}
Notice that the random variables $L_m$, $m \geq 0$, as well as the variables $K_m$, $m \geq 0$, are not 
independent.

\medskip

\textbf{Example 5.} Let us calculate $\mathbb{P}\left\{L_1 = 0\right\}$.
In view of \eqref{J5}, \eqref{J8}, and \eqref{J9}
(as well as \eqref{A1} and the independence of $X_0^1$ and $X^1_1$)  we have
\begin{align}
&\mathbb{P}\left\{L_1 = 0\right\} = \mathbb{P}\left\{X^{L_0}_1 = 0\right\}
=\sum_{k=1}^{\infty} \mathbb{P}\left\{X^{L_0}_1 = 0\,|\, L_0 = k\right\} \mathbb{P}\left\{L_0 = k\right\}
\nonumber
\\
& \quad = \mathbb{P}\left\{X^1_1 = 0\,|\, X_0^1 = 0\right\} (1-\beta)
+ (1-\beta)\sum_{k=2}^{\infty} \mathbb{P}\left\{X^k_1 = 0\,|\, L_0 = k\right\} \beta^{k-1}
\nonumber
\\
& \quad = (1-\alpha) (1-\beta)
+ (1-\beta)\sum_{k=2}^{\infty} \mathbb{P}\left\{X^k_1 = 0\,|\, L_0 = k\right\} \beta^{k-1}.
\label{J11}
\end{align}
Now, in view of \eqref{J6} we have that
\begin{equation}
\mathbb{P}\left\{X^k_1 = 0\,|\, L_0 = k\right\}
= \mathbb{P}\left\{X^k_1 = 0\,|\, X_0^1 = X_0^2 = \cdots = X_0^{k-1} = 1, X^k_0 = 0\right\}.
\label{J12}
\end{equation}
Proceeding as in Example 4 we get
\begin{align}
&\mathbb{P}\left\{X^k_1 = 0\,|\, X_0^1 = X_0^2 = \cdots = X_0^{k-1} = 1, X^k_0 = 0\right\}
\nonumber
\\
& \quad = (1-\alpha) w_{k-1} + (1-\gamma) (1 - w_{k-1}) = (\gamma - \alpha) w_{k-1} + 1-\gamma,
\label{J13}
\end{align}
where $w_k$ is given by \eqref{B16}.

Finally, by using \eqref{J13}, \eqref{J12}, and \eqref{B16} in \eqref{J11} we obtain (by summing the 
resulting geometric series)
\begin{equation}
\mathbb{P}\left\{L_1 = 0\right\}
= \frac{1 - (\gamma - \alpha) \beta - \alpha}{1 - (\gamma - \alpha) \beta}.
\label{J14}
\end{equation}
\hfill $\diamondsuit$



%
%

\medskip

\section{Epilogue} The model presented in this paper can be extended/enriched in various ways.
For instance:

(i) We can assume that the fire starts at many points, say at the points
$(j, -1)$, $k \in F_x$ and $(-1,k)$,
$k \in F_y$, where $F_x, F_y \subset \mathbb{N}$. In this case, the ``boundary conditions" for 
$X(j, k)$ are
\begin{equation*}
X(j, -1) = X(-1, k) = 1,
\quad
j \in F_x, \ k \in F_y,
\end{equation*}
and
\begin{equation*}
X(j, -1) = X(-1, k) = 0,
\quad
j \in \mathbb{N} \smallsetminus F_x, \ k \in \mathbb{N} \smallsetminus F_y.
\end{equation*}
With such boundary conditions we can furthermore assume that $X(-1, -1) = 0$ and, 
also, $X(j, k) = 0$ whenever $j \leq -2$ or $k \leq -2$.

(ii) The status $X(j, k)$ of the tree at the point $(j, k)$ can take more than two values, say
$X(j, k) \in \{0, 1, \ldots, m\}$ to take into account the case where the tree at $(j, k)$ is partially burnt (i.e. burnt to a certain degree).

(iii) We can consider a higher dimensional model. For example we can consider a corrosion model:
Each point $(j, k, \ell) \in \mathbb{N}^3$ is associated to a molecule (or atom, or ion) and the status $X(j, k, \ell)$ of the molecule at $(j, k, \ell)$ is affected only by 
the status of the neighboring molecules at $(j-1, k, \ell)$, $(j, k-1, \ell)$, and $(j, k, \ell-1)$.

(iv) Finally, one might investigate whether there is a continuous analog of the 
presented model.

\section{APPENDIX: An one-dimensional model}
Let us discuss briefly the (almost trivial) one-dimensional case, where the forest is represented by 
the set $\mathbb{N} = \{0, 1, 2, \ldots\}$ and each point $j \in \mathbb{N}$ is associated to a 
tree. The status of the tree at $j$ is $X(j)$, so that if $X(j) = 1$, the tree at $j$ is burnt, 
while if $X(j) = 0$, then the tree at $j$ is not burnt.

Here the model goes as follows: The status of the tree at $j$ is affected only by the status of the neighboring tree
at $j-1$ so that
\begin{equation}
\mathbb{P}\{X(j) = 1\,|\, X(j-1) = 1\} = p
\qquad
\text{and}
\qquad
\mathbb{P}\{X(j) = 1\,|\, X(j-1) = 0\} = 0,
\label{D1}
\end{equation}
where
\begin{equation*}
0 \leq p \leq 1.
\end{equation*}
The parameter (probability) $p$ depends on the strength of the wind. The smaller the 
magnitude of the velocity of the wind, the closer to $0$ is $p$.

We assume that the fire starts at the point $0$ with a given probability $r \in [0, 1]$. Thus the boundary condition for $X(j)$ is
\begin{equation}
\mathbb{P}\{X(0) = 1\} = r
\label{D3}
\end{equation}
(this is slightly more general than the condition $X(0) = 1$).

Notice that our assumptions imply
\begin{equation}
\{X(j) = 1\} \subset \{X(j-1) = 1\},
\label{D6}
\end{equation}
hence
\begin{equation}
\{X(j) = 1\} = \{X(0) = 1\} \cap \{X(1) = 1\} \cap \cdots \cap \{X(j) = 1\}.
\label{D7}
\end{equation}

From \eqref{D1} and \eqref{D3} we get that
\begin{equation*}
\mathbb{P}\{X(1) = 1\} = \mathbb{P}\{X(1) = 1\,|\, X(0) = 1\} \mathbb{P}\{X(0) = 1\} = rp,
\end{equation*}
\begin{equation*}
\mathbb{P}\{X(2) = 1\} = \mathbb{P}\{X(2) = 1\,|\, X(1) = 1\} \mathbb{P}\{X(1) = 1\}
 = rp^2,
\end{equation*}
and, in general,
\begin{equation}
\mathbb{P}\{X(j) = 1\} = rp^j,
\qquad
j \geq 0.
\label{D5}
\end{equation}

Let
\begin{equation}
L = \inf\left\{n \geq 0 : X(n) = 0\right\},
\label{D8}
\end{equation}
so that
\begin{equation}
X(n+L) = 0
\qquad
\text{for all }
n \geq 0. 
\label{D9}
\end{equation}
In fact, $L$ is the (total) number of burnt trees.

From \eqref{D3} and \eqref{D5} we get
\begin{equation}
\mathbb{P}\{L = 0\} = \mathbb{P}\{X(0) = 0\} = 1-r
\label{D10a}
\end{equation}
and
\begin{align}
\mathbb{P}\{L = k\} &= \mathbb{P}\{X(k) = 0, X(k-1) = 1\}
\nonumber
\\
&= \mathbb{P}\{X(k) = 0\,|\, X(k-1) = 1\} \mathbb{P}\{X(k-1) = 1\}
\nonumber
\\
&= r (1-p) p^{k-1},
\qquad
k \geq 1.
\label{D10}
\end{align}
Notice that the conditional distribution of the random variable $L$ given
$\{X(0) = 1\}$ is geometric with parameter 
$1-p$. Thus, unless $p = 1$ in which case the whole forest will be burnt,
we have that $L < \infty$ almost surely (which means that most of the forest will be 
spared).

Finally, by a direct calculation, or by conditioning on $\{X(0) = 0\}$ and $\{X(0) = 1\}$ and
by invoking the properties of the geometric distribution we easily obtain that
\begin{equation}
\mathbb{E}[L] = \frac{r}{1-p}
\qquad
\text{and}
\qquad
\mathbb{V}[L] = \frac{r(1+p-r)}{(1-p)^2}.
\label{D12}
\end{equation}


%
%
%
%
%
%
%
%
%
%

\end{document}